# Bayesian Renewables Scenario Generation via Deep Generative Networks

Yize Chen, Pan Li, and Baosen Zhang
Department of Electrical Engineering
University of Washington
Seattle, Washington 98195
Email: {yizechen, pli69, zhangbao}@uw.edu

*Abstract*—We present a method to generate renewable scenarios using Bayesian probabilities by implementing the Bayesian generative adversarial network (Bayesian GAN), which is a variant of generative adversarial networks based on two interconnected deep neural networks. By using a Bayesian formulation, generators can be constructed and trained to produce scenarios that capture different salient modes in the data, allowing for better diversity and more accurate representation of the underlying physical process. Compared to conventional statistical models that are often hard to scale or sample from, this method is model-free and can generate samples extremely efficiently. For validation, we use wind and solar times-series data from NREL integration data sets to train the Bayesian GAN. We demonstrate that proposed method is able to generate clusters of wind scenarios with different variance and mean value, and is able to distinguish and generate wind and solar scenarios simultaneously even if the historical data are intentionally mixed.

## I. INTRODUCTION

The stochastic and intermittent nature of renewable resources has brought new challenges to the scheduling, operation, and planning of power systems. One popular framework to capture these uncertainties in renewable generation is to use a set of scenarios, each representing a possible time series realization of the random physical process [1], [2]. These scenarios can then be used in a variety of optimization problems, including stochastic economic dispatch, unit commitment, operation of wind farms, storage management, trading and planning (see, e.g. [3], [4], [5], [6] and the references within).

A key requirement of the generated scenarios is that they accurately reflect the spatial/temporal patterns in the physical generation process and exhibit enough diversity to capture a wide range of behaviors. For instance, most land-based wind farms have diurnal patterns but may also output very little power over a long period of time, and a set of scenarios should ideally capture both phenomena. To capture underlying stochastic processes, several model-based approaches have been proposed, where the process is assumed to have some parametrized probability distribution $\mathscr{P}$, and historical data is used to learn these parameters [7], [8]. For example, in [2], [9], copula methods are first applied to model the distribution and correlation of power generation time-series, then scenarios are generated via a sampling method. In [1], an autoregressive moving average (ARMA) model is first learned and then used to generate spatiotemporal scenarios using power generation profiles at each renewables generation site. However, the time-varying and nonlinear dynamics of weather along with the complex spatial and temporal interactions make model-based approaches difficult to apply and hard to scale, especially when multiple renewable power plants are considered. In particular, even if the exact distribution of the stochastic process is known, sampling from it is usually nontrivial and time consuming. Moreover, some of previous methods depend on certain probabilistic forecasts as inputs, which could limit the diversity of the generated scenarios and under-explore the overall variability of renewable resources [10].

In [11], we proposed a deep generative machine learning framework to capture and learn the power generation dynamics. We adopted the Generative Adversarial Networks (GAN), which is a *data-driven, model-free* approach that generates new scenarios directly from historical data without explicitly fitting a probabilistic model. This approach is easily scalable, and outperforms existing state-of-the-art model-based approaches especially in the setting of multiple correlated renewable generators. However, despite its success, if the historical data contains several distinct modes (e.g., high wind and low wind profiles), the generated scenarios can contain a mixture of these modes. These "mixed" scenarios maybe inappropriate for the subsequent optimization problems and need to be filtered out with an additional post-processing step.

To overcome this challenge, we extend the results of [11] by introducing a Bayesian formulation into the GAN training process. This Bayesian formulation allows us to use the GAN architecture to directly find a *set of generators* to fully capture different modes in the historical data. This approach is based on the result in [12] which introduced the Bayesian GAN. During the training process, we sample neural network weights from a prior distribution to form a group of particular generative models. By exploring the posterior distribution over the parameters of both the generator and discriminator, we are able to obtain this set of generators. The framework of proposed method for scenario generation is plotted in Fig. 1.

The following sections are organized as follows: in Section. II we describe the preliminaries on problem formulation as well as the motivation and Bayesian inference idea for incorporating Bayesian GAN model; in Section. III we introduce our generative model setup and training algorithm; numerical simulations are performed and evaluated in Section. IV. We will show each trained generator is able to generate specific scenarios under certain dynamics for several scenario generation tasks.

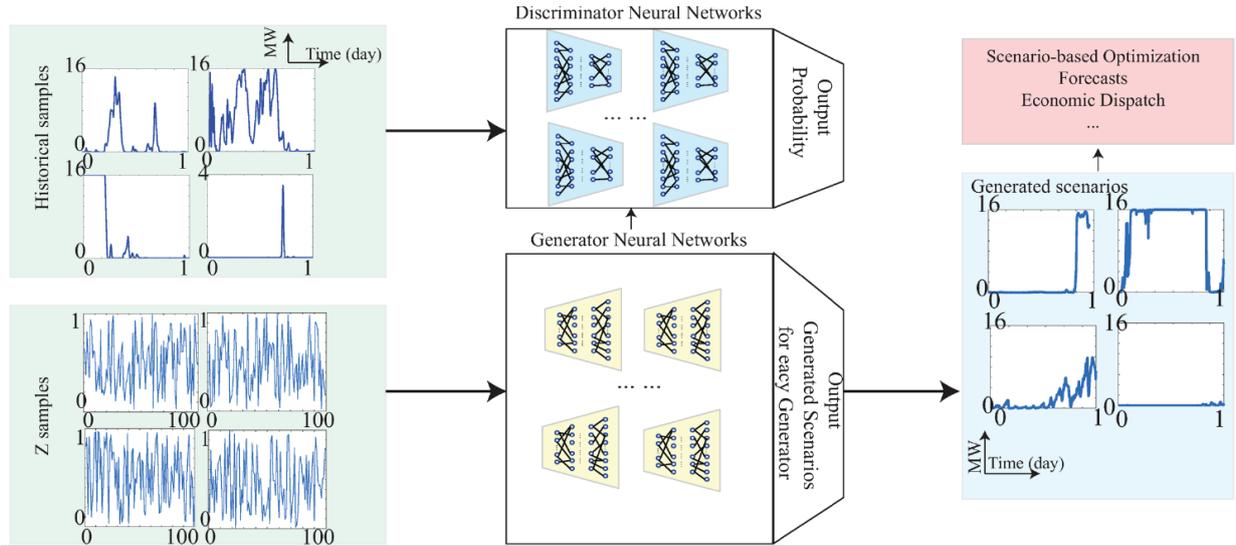

Fig. 1. Deep generative model framework for our proposed scenario generation method using Bayesian information. By using sampling method to explore the full posterior over the discriminator and generator neural networks' weights, we are able to use a set of generators producing scenarios of different modes that depict different dynamic behaviors for renewables generation process.

## II. Preliminaries

In this section we present our problem formulation to generate scenarios, and discuss some arising issues which motivate us to enhance its performance by using a Bayesian formulation.

### A. Problem Formulation

Consider a set of historical data for a group of renewable resources at $N$ sites. For site $j$, let $\mathbf{x}_j$ be the vector of historical data indexed by time, $t = 1, \ldots, T$, and $j$ ranges from 1 to $N$.

Our objective is to train a generative model based on GANs (with neural network parameters $\theta$) by utilizing historical power generation data as the training set. And we are interested in two scenario generation problem:

*Single Site Scenario Generation:* In this case, we want to generate a group of scenarios representing the dynamics of certain renewables generation farm using $\mathbf{x}_j$.

*Spatiotemporal Scenario Generation:* In this case, we want to generate spatiotemporal scenarios for a group of renewables generation sites $\{\mathbf{x}_j\}, j = 1, \ldots, N$.

Generated scenarios should be capable of describing the same stochastic processes as training samples and exhibiting a variety of different modes representing all possible variations and patterns seen during training.

### B. Drawbacks of a Single Generator

The generated scenarios should capture the inherent diversity (both spatial and temporal) of renewables generation. However, the following two problems can potentially occur for a generator trained with GAN or in more conventional model-based methods where a single model is trained:

1) Generated scenarios do not capture the diversity in renewable generation scenarios, i.e., most of the generated scenarios resemble each other and do not provide enough information for possible range realizations of power generation process.
2) It is hard to separate out distinct behaviors through a single generator. Since the generator finds the mapping from input noisy distribution to overall joint distribution of the power output, there is no straightforward way to look into scenarios of specific characteristics (e.g., larger power generation values or frequent ramp events).

Both problems are the result of neural networks learning from the most frequent training samples and not generalizing. These problems can mislead decisions in power system because the generated scenarios are not able to depict the true distribution of the renewable power output. The following sections describe a Bayesian method to place a prior distribution on the neural network weights. Moreover, by conducting appropriate sampling techniques, the proposed method is able to track multiple local optima, thus enhance the capability to capture mode diversity in targeted data distribution. In Section II-C, we briefly review the basics of Bayesian inference and leave more details into Section III.

### C. Bayesian inference

Given a neural network $\Omega$, with data samples $\{\mathbf{x}_j\}, j = 1, \ldots, N$ and learnable parameter $\theta$ (e.g., neural network weights), the traditional way to find optimal values for $\theta$ is to maximize the data likelihood:

$$\theta^* = \arg\max_{\theta} \prod_{j=1}^{N} p_{\Omega}\{x_j | \theta\}, \tag{1}$$

where $p(\cdot)_\Omega$ denotes likelihood of data $x_j$ that depends on the structure of the neural network.

By adding a prior on $\theta$, the objective for the neural network becomes:

$$\theta^* = \arg\max_\theta \prod_{j=1}^{N} p_\Omega\{x_j|\theta\} p\{\theta|\gamma\}, \quad (2)$$

where $p\{\theta|\gamma\}$ is the prior distribution on $\theta$ with parameter $\gamma$.

In Section III, we discuss more details on the structure of a traditional GAN and how we incorporate Bayesian information to improve the performance of GAN.

## III. BAYESIAN GAN

In this section we proceed with more details on GAN setup, and present a variant of GAN introduced in [12] that uses Bayesian inference to improve the generative performance of GAN by training a group of neural networks together.

### A. GANs with Wasserstein Distance

The general architecture of proposed method we use is shown in Fig. 1. Assume observations $x_j^t$ for times $t \in T$ of renewable power are available for each power plant $j$, $j = 1,...,N$. Let the true distribution of the observation be denoted by $\mathbb{P}_X$. We have access to a group of noise vector input $z$ under a known distribution $Z \sim \mathbb{P}_Z$ which can be easily sampled from (e.g., from jointly Gaussian). Our goal is to transform ple $z$ drawn from $\mathbb{P}_Z$ such that it follows $\mathbb{P}_X$ (without ever learning $\mathbb{P}_X$ explicitly). This is accomplished by simultaneously training two deep neural networks: the generator network and the discriminator network. Let $G$ denote the generator function parametrized by $\theta^{(G)}$, which we write as $G(\cdot;\theta^{(G)})$; and let $D$ denote the generator function parametrized by $\theta^{(D)}$ which we write as $D(\cdot;\theta^{(D)})$. Here, $\theta^{(G)}$ and $\theta^{(D)}$ are the weights of two neural networks, respectively. For convenience, we sometimes suppress the symbol $\theta$.

*Generator:* During the training process, the generator is trained to take a batch of inputs and by taking a series of up-sampling operations by neurons of different functions to output realistic scenarios. Suppose that $Z$ is a random variable with distribution $\mathbb{P}_Z$. Then $G(Z;\theta^{(G)})$ is a new random variable, whose distribution is denoted as $\mathbb{P}_G$.

*Discriminator*: The discriminator is trained simultaneously with the generator. It takes input samples either coming from real historical data or coming from generator, and by taking a series of operations of down-sampling using another deep neural network, it outputs a continuous value $p_{real}$ that measures to what extent the input samples belong to $\mathbb{P}_X$:

$$p_{real} = D(x;\theta^{(D)}) \quad (3)$$

where $x$ may come from $\mathbb{P}_{data}$ or $\mathbb{P}_Z$. The discriminator is trained to learn to distinguish between $\mathbb{P}_X$ from $\mathbb{P}_G$, and thus to maximize the difference between $\mathbb{E}[D(X)]$ (real data) and $\mathbb{E}[D(G(Z))]$ (generated data).

With the objectives for discriminator and generator defined, we need to formulate loss function $L_G$ for generator and $L_D$ for discriminator to train them (i.e., update neural networks' weights based on the defined loss function). Following this guideline and the loss defined in [13], we can write $L_D$ and $L_G$ as followed:

$$L_G = -\mathbb{E}_Z[D(G(Z))] \quad (4a)$$
$$L_D = -\mathbb{E}_X[D(X)] + \mathbb{E}_Z[D(G(Z))]. \quad (4b)$$

Since a large discriminator output means the sample is more realistic, the generator will try to minimize the expectation of $-D(G(\cdot))$ by varying $G$ (for a given $D$), resulting in the loss function in (4a). On the other hand, for a given $G$, the discriminator wants to minimize the expectation of $D(G(\cdot))$, and the same time maximizing the score of real historical data. This gives the loss function in (4b).

We then combine (4a) and (4b) to form a two-player iterative minimax game with the value function $V(G,D)$:

$$\min_{\theta^{(G)}} \max_{\theta^{(D)}} V(G,D) = \mathbb{E}_X[D(X)] - \mathbb{E}_Z[D(G(Z))] \quad (5)$$

where $V(G,D)$ is the negative of $L_D$.

For more training details, we refer readers to [14], [13], [11]. As training moves on and goes near to the optimal solution, $G$ is able to generate samples that look as realistic as real data with a small $L_G$ value, while $D$ is unable to distinguish $G(z)$ from $\mathbb{P}_X$ with large $L_D$. Eventually, we are able to learn an unsupervised representation of the probability distribution of renewables scenarios from the output of $G$.

More formally, the minimax objective (5) of the game can be interpreted as the dual of the so-called Wasserstein distance (Earth-Mover distance) [15], [16]. Under the setting of 4a 4b, we are precisely trying to get two random variables, $\mathbb{P}_X(D(X))$ and $\mathbb{P}_Z(D(G(X)))$, to be close to each other. It turns out that

$$W(D(X),D(G(Z))) = \sup_{\theta^{(D)}}\{\mathbb{E}_X[D(X)] - \mathbb{E}_Z[D(G(Z))], \quad (6)$$

which is a natural connection to (5). The expectations can be computed as empirical means in the mini-batch updates for neural network training.

### B. Bayesian GAN

As motivated by Section II-C, we want to use Bayesian information to enhance the performance of GAN and avoid mode collapse in generated scenarios. Moreover, we would like to examine that by training potentially a group of generators together, each individual generator is learning unique dynamics of power generation process. Inspired by [12], we interpret losses defined in (4a)(4b) as the negative log likelihood for $\theta^{(G)}$ and $\theta^{(D)}$, then we formulate the posteriors of $\theta^{(G)}$ and $\theta^{(D)}$ similarly to (2):

$$\log p\{\theta^{(G)}|\theta^{(D)}\} = \mathbb{E}_Z[D(G(Z,\theta^{(G)}),\theta^{(D)})] + \log p\{\theta^{(G)}|\gamma^{(G)}\}, \quad (7)$$

and

$$\log p\{\theta^{(D)}|\theta^{(G)}\} = (\mathbb{E}_X[D(X,\theta^{(D)})] \\ - \mathbb{E}_Z[D(G(Z,\theta^{(G)}),\theta^{(D)})]) + \log p\{\theta^{(D)}|\gamma^{(D)}\}, \quad (8)$$

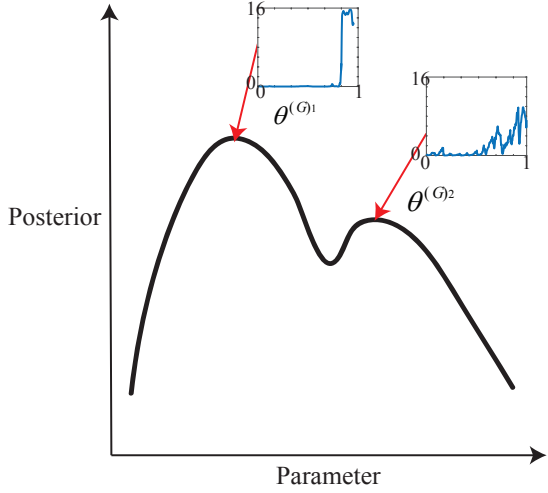

Fig. 2. The sketch of the multimodal posterior over the generator neural networks' weights. By using two different generators, they are able to generate scenarios of unique behaviors.

where $p\{\theta^{(G)}|\gamma^{(G)}\}$ and $p\{\theta^{(D)}|\gamma^{(D)}\}$ are priors over $\theta^{(G)}$ and $\theta^{(D)}$. The priors over $\theta^{(G)}$ and $\theta^{(D)}$ represent the initial estimation of $\theta^{(G)}$ and $\theta^{(D)}$ with hyperparameters $\gamma^{(G)}$ and $\gamma^{(D)}$. Once trained with large number of historical samples **x**, the constructed posteriors of $\theta^{(G)}$ and $\theta^{(D)}$ represent the refined estimation of $\theta^{(G)}$ and $\theta^{(D)}$ based on samples.

We can obtain a valid generator by maximizing the posterior distribution. However, that gives us one particular generator (and one discriminator), which may not be able to capture all the dynamics of the power generation. To obtain a group of generators that reproduce the different dynamics, we conduct iterative posterior sampling instead of posterior maximization:

$$\theta^{(D)}|\theta^{(G)} \sim p\{\theta^{(D)}|\theta^{(G)}\}, \theta^{(G)}|\theta^{(D)} \sim p\{\theta^{(G)}|\theta^{(D)}\}. \quad (9)$$

There are several reasons for adopting posterior sampling:

- The imposed prior on $\theta^{(G)}$ and $\theta^{(D)}$ can be regarded as the regularization term, which reduces overfitting on training data.
- By iterative sampling over posterior distribution, we are able to capture multi-mode optima in posterior distribution and further avoid mode collapse, as shown in Fig. 2. A traditional GAN can be seen as iteratively maximizing the posteriors over $\theta^{(G)}$ and $\theta^{(D)}$ with a uniform prior.

### C. Algorithm

To sample from posterior distributions in (7) and (8), we adopt the Stochastic Gradient Hamiltonian Monte Carlo (SGHMC) in [12] for posterior sampling. The whole algorithm for training Bayesian GAN is presented in Algorithm 1, while the code is public on Github [1].

[1] https://github.com/chennnnnyize/BayesianRenewablesGAN

**Algorithm 1** Training Bayesian GAN using SGHMC, adopted from Algorithm 1 in [12].

**Require:** Learning rate $\alpha$, friction term $\eta$, clipping parameter $c$, batch size $m$, Number of iterations for discriminator per generator iteration $n_{discri}$, MC iterations for discriminator $N_d$ and for generator $N_g$, number of SGHMC samples $M$

**Require:** Initial weights $\theta^{(D)}$ for discriminator and $\theta^{(G)}$ for generator, initial training samples $N = 0$

**while** $\theta^{(D)}$ has not converged **do**
  **for** number of MC iterations $N_d$ **do**
    **for** $t = 0, ..., n_{discri}$ **do**
      *# Update parameter for Discriminator*
      Sample mini-batch data:
      $\{(x^{(i)}, y^{(i)})\}_{i=1}^m from \mathbb{P}_X$ $\{(z^{(i)}, y^{(i)})\}_{i=1}^m from \mathbb{P}_Z$
      Update discriminator nets by running SGHMC for $M$ times:
      $N = N + m$
      $g_{\theta^{(D)}} \leftarrow \nabla_{\theta^{(D)}}[-\frac{1}{m}\sum_{i=1}^m D(x^{(i)}|y^{(i)}) + \frac{1}{m}\sum_{i=1}^m D(G(z^{(i)}|y^{(i)})) + \frac{1}{N}\log p(\theta^{(D)}|\gamma^{(D)}) + \mathbf{n}], \mathbf{n} \sim \mathcal{N}(0, 2\eta\alpha I)$.
      $\theta^{(D)} \leftarrow \theta^{(D)} - \alpha \cdot RMSProp(\theta^{(D)}, g_{\theta^{(D)}})$
      $\theta^{(D)} \leftarrow clip(w, -c, c)$
    **end for**
  **end for**
  **for** number of MC iterations $N_g$ **do**
    *# Update parameter for Generator*
    Update generator nets by running SGHMC for $M$ times:
    $N = N + m$
    $g_{\theta^{(G)}} \leftarrow \nabla_{\theta^{(G)}}[\frac{1}{m}\sum_{i=1}^m D(G(z^{(i)}|y^{(i)})) + \frac{1}{N}\log p(\theta^{(G)}|\gamma^{(G)}) + \mathbf{n}], \mathbf{n} \sim \mathcal{N}(0, 2\eta\alpha I)$.
    $\theta^{(G)} \leftarrow \theta^{(G)} - \alpha \cdot RMSProp(\theta^{(G)}, g_{\theta^{(G)}})$
  **end for**
**end while**

## IV. SIMULATION RESULTS

In order to test the performance of our proposed framework for scenario generation, we set up our case studies based on wind power and PV generation data published by the National Renewable Energy Lab [2]. The resolution of the data is 5 minutes. We consider the following settings and test if Bayesian GAN is able to produce diverse, particular scenarios under the following circumstances:

1) Training input is composed of a mixture of wind power and PV samples. We want to examine whether we can train two generators such that one generates scenarios of wind power and the other generates scenarios of PV using (7) and (8). This serves as a sanity check for our proposed method.
2) Training data is composed of group of spatiotemporal wind samples, i.e., group of wind farms under different spatiotemporal correlation. We want to check

[2] https://www.nrel.gov/grid

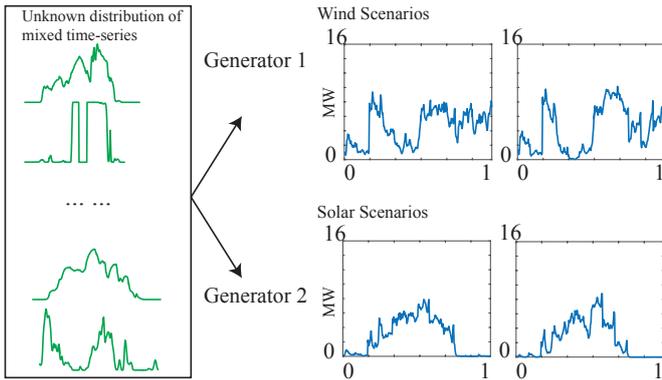

Fig. 3. By training two generators together ($J = 2, M = 2$ specifically) using Bayesian GAN, 1-day scenarios generated by each generator capture the distinct wind and solar power generation dynamics respectively.

   if Bayesian inference would help distinguish specific correlations.
3) Training input only includes wind samples. We would like to evaluate if different trained generators could generate scenarios with unique dynamics (e.g., ramp events, power generation mean value).

### A. Co-generation of Wind and PV Scenarios

We first test the performance of Bayesian GAN by artificially blending wind and PV scenarios into the training data. The input data is 7-year collection of PV generation profiles and wind power generation profiles. Both datasets are split into 24-hour samples and all samples are shuffled to form an overall mixed distribution of solar and wind generation data.

This is a toy model, since we want to examine when trained properly, different generators in Bayesian GAN are learning the different modes of the input distribution. An interesting observation is that the two generators generate the respective wind and PV scenarios efficiently as shown in Fig 3. In comparison, if we train the normal GAN model, the generator finds the map from input noise to the single centered mode of the input distribution, thus generated scenarios are a mixture of solar and wind power generation samples.

### B. Spatiotemporal Scenario Generation

In this subsection, we examine the model performance on generating spatiotemporal scenarios. We extend the wind datasets we used in Section. IV-A, and add a group of 24 wind farms' historical observations. These 24 wind farms are located in the State of North Dakota, US, and exhibit quite different spatiotemporal correlations compared to previous group of wind farms located in Texas. Instead of inputting single farm power generation time-series, we input the 24 wind farms' 24-hour scenarios as a matrix into GAN. Similar to the shuffling step we take in Section. IV-A, we train 2 generators using $M = 2, J = 2$ to check if generated scenarios by two distinct generators exhibit distinct spatiotemporal correlations.

To evaluate if generated spatiotemporal scenarios observe the similar correlations, we calculate the Pearson correlation coefficient $\rho_{i,j}$ for wind farm $\mathbf{x}_i$ and $\mathbf{x}_j$ [17]. The color map for two groups' historical data and generated scenarios' correlation coefficient matrix are plotted in Fig 4. We could observe similar correlation coefficient matrix structure for each pair of generated scenarios and original training data. Yet due to the similarity of input samples (e.g., power generation peak values, mean values, ramp events), the two generators can not totally distinguish these two processes. It still provides us a useful tool for extracting spatiotemporal patterns existing in renewables generation profiles.

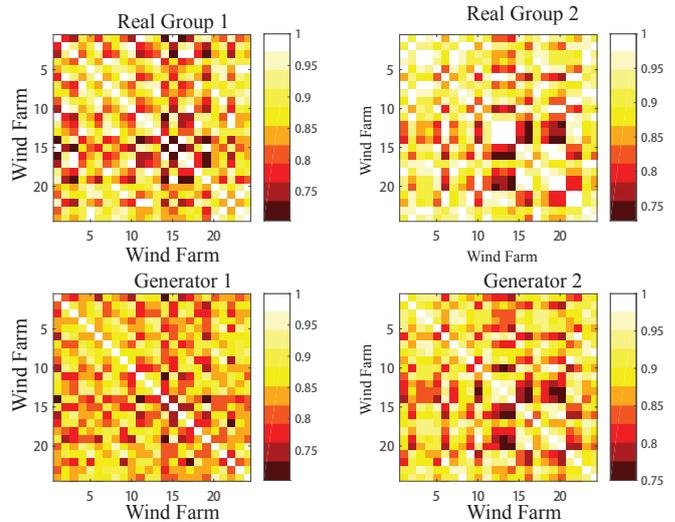

Fig. 4. The spatiotemporal scenarios for a group of 24 wind farms located in State of Texas and North Dakota respectively. By using Bayesian GAN, Generator 1 learns Real data Group 1, Generator 2 learns Real data Group 2 spatiotemporal correlation respectively.

### C. Single-Farm Wind Scenario Generation

In this group of simulations we test whether Bayesian GAN is able to train distinct generators that represent unique wind power generation dynamics.

We collect 7-year of wind generation data for nearby 24 wind farms in State of Texas, US, and form them as a group to represent the historical power generation dynamics for that location. All of the wind farms have a capacity of $16MW$. Once these samples are fed into our generative model framework, we train Bayesian GAN using $J = 4, M = 4$ till the discriminator loss converges.

In Fig. 5 we illustrate the simulation result using these 4 generators. In Fig. 5(a) and Fig. 5(b) we could observe that scenarios generated by generator 2 have much smaller variance and mean values compared with the other 3 generators. Thus this trained generator could be used for generating scenarios for mild days. On the contrary, scenarios generated by generator 1 represent most of the scenarios with frequent ramp events and high generation output. Such differences of generated scenarios' dynamics are also depicted in the randomly selected output scenarios in Fig. 5(c).

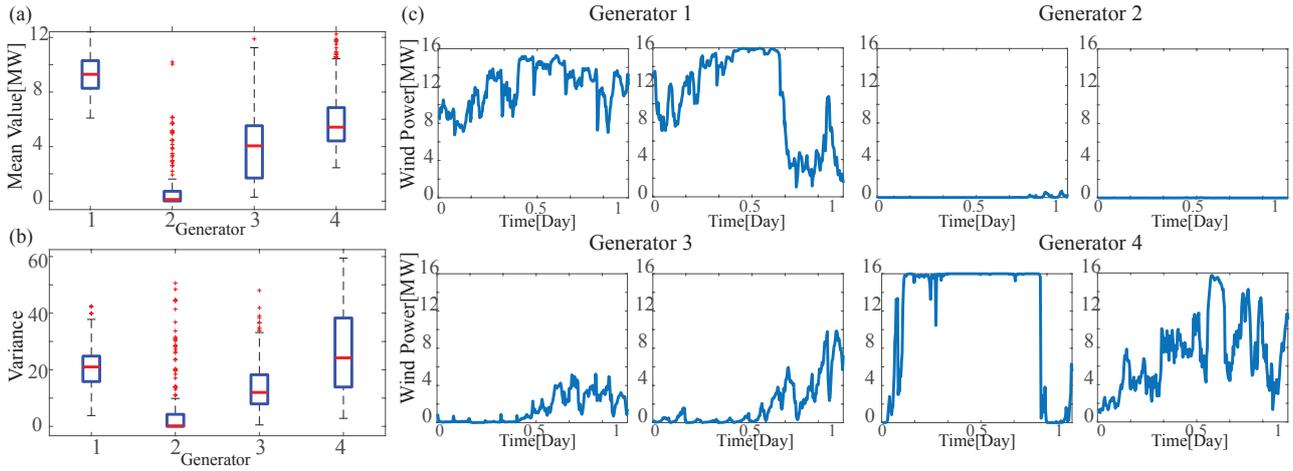

Fig. 5. Scenarios generated with $J = 4, M = 4$. (a) and (b) are the boxplots for generated scenarios' mean value (MW) and variance for four different generators respectively; in (c) we randomly select and plot the distinct generated scenarios for four generators.

## V. CONCLUSION AND DISCUSSION

Scenario generation is an essential tool for decision-making in power grids with high penetration of renewables. In this paper, we incorporate Bayesian information into the Generative Adversarial Networks and implement the Bayesian GAN to capture inherent diversity in clusters of scenarios. Our proposed method leverages the power of deep neural networks and large sets of historical data to perform the task for directly generating scenarios conforming to the same distribution of historical data, without explicitly modeling of the distribution. The simulation results based on NREL wind and solar data show that proposed method is able to capture the inherent diversity in the data distribution and reproduce each of them.

Our proposed method can be applied to most stochastic processes of interest in power systems. In addition, as the method uses a feed forward neural network structure, it does not require sampling and can be scaled easily to systems with a large number of uncertainties.

## REFERENCES

[1] J. M. Morales, R. Minguez, and A. J. Conejo, "A methodology to generate statistically dependent wind speed scenarios," *Applied Energy*, vol. 87, no. 3, pp. 843–855, 2010.
[2] P. Pinson, H. Madsen, H. A. Nielsen, G. Papaefthymiou, and B. Klöckl, "From probabilistic forecasts to statistical scenarios of short-term wind power production," *Wind energy*, vol. 12, no. 1, pp. 51–62, 2009.
[3] Q. P. Zheng, J. Wang, and A. L. Liu, "Stochastic optimization for unit commitmentâĂŤa review," *IEEE Transactions on Power Systems*, vol. 30, no. 4, pp. 1913–1924, 2015.
[4] H. Park and R. Baldick, "Transmission planning under uncertainties of wind and load: Sequential approximation approach," *IEEE Transactions on Power Systems*, vol. 28, no. 3, pp. 2395–2402, 2013.
[5] A. Papavasiliou, S. S. Oren, and R. P. O'Neill, "Reserve requirements for wind power integration: A scenario-based stochastic programming framework," *IEEE Transactions on Power Systems*, vol. 26, no. 4, pp. 2197–2206, 2011.
[6] H. Ding, Z. Hu, and Y. Song, "Stochastic optimization of the daily operation of wind farm and pumped-hydro-storage plant," *Renewable Energy*, vol. 48, pp. 571–578, 2012.
[7] J. Dupačová, G. Consigli, and S. W. Wallace, "Scenarios for multistage stochastic programs," *Annals of operations research*, vol. 100, no. 1-4, pp. 25–53, 2000.
[8] D. Lee and R. Baldick, "Load and wind power scenario generation through the generalized dynamic factor model," *IEEE Transactions on Power Systems*, vol. 32, no. 1, pp. 400–410, 2017.
[9] E. K. Hart and M. Z. Jacobson, "A monte carlo approach to generator portfolio planning and carbon emissions assessments of systems with large penetrations of variable renewables," *Renewable Energy*, vol. 36, no. 8, pp. 2278–2286, 2011.
[10] G. Sideratos and N. D. Hatziargyriou, "Probabilistic wind power forecasting using radial basis function neural networks," *IEEE Transactions on Power Systems*, vol. 27, no. 4, pp. 1788–1796, 2012.
[11] Y. Chen, Y. Wang, D. S. Kirschen, and B. Zhang, "Model-free renewable scenario generation using generative adversarial networks," *IEEE Transactions on Power Systems*, 2018.
[12] Y. Saatchi and A. Wilson, "Bayesian GAN," in *Advances in Neural Information Processing Systems*, 2017, pp. 3625–3634.
[13] M. Arjovsky, S. Chintala, and L. Bottou, "Wasserstein GAN," *arXiv preprint arXiv:1701.07875*, 2017.
[14] I. Goodfellow, J. Pouget-Abadie, M. Mirza, B. Xu, D. Warde-Farley, S. Ozair, A. Courville, and Y. Bengio, "Generative adversarial nets," in *Advances in neural information processing systems*, 2014, pp. 2672–2680.
[15] C. Villani, *Optimal transport: old and new*. Springer Science & Business Media, 2008, vol. 338.
[16] Y. Rubner, C. Tomasi, and L. J. Guibas, "The earth mover's distance as a metric for image retrieval," *International journal of computer vision*, vol. 40, no. 2, pp. 99–121, 2000.
[17] P. Pinson and R. Girard, "Evaluating the quality of scenarios of short-term wind power generation," *Applied Energy*, vol. 96, pp. 12–20, 2012.